\documentclass{article}[11pt]

\usepackage[ruled,vlined,linesnumbered]{algorithm2e}
\usepackage{amsfonts,amsmath,tikz}

\usepackage{graphicx}
\usepackage{amsmath,amsfonts,amssymb,latexsym}
\usepackage{enumerate}
\usepackage{setspace}
\usepackage{color}
\usepackage[cp1250]{inputenc}
\usepackage[width=145mm,height=205mm,left=35mm,foot=10mm]{geometry}

\newtheorem{thm}{Theorem}[section]
\newtheorem{prp}[thm]{Proposition}

\newtheorem{cor}[thm]{Corollary}
\newtheorem{rem}[thm]{Remark}

\newcommand{\qed}{\hfill ~$\square$\bigskip}
\newcommand{\proof}{\noindent{\bf Proof.} }

\newcommand{\cp}{\,\square\,}
\newcommand{\vertex}{\node[vertex]}
\tikzstyle{vertex}=[circle, draw, inner sep=0pt, minimum size=6pt]

\newcommand{\ggrz}{\gamma_{gr}^{\textrm{Z}}}
\newcommand{\ggr}{\gamma_{gr}}
\newcommand{\ggrt}{\gamma_{gr}^{t}}
\newcommand{\ggrl}{\gamma_{gr}^{L}}

\newcommand{\cH}{{\cal H}}

\begin{document}

\title{Grundy dominating sequences and zero forcing sets}

\author{
Bo\v stjan Bre\v sar $^{a,b}$ \and Csilla Bujt{\'a}s $^{c,d}$ \and
Tanja Gologranc $^{a,b}$ \and Sandi Klav\v zar $^{e,a,b}$ \and Ga\v sper
Ko\v smrlj $^{b,f}$ \and Bal{\'a}zs Patk{\'o}s $^{d}$ \and Zsolt Tuza
$^{c,d}$ \and M{\'a}t{\'e} Vizer $^{d}$ }

\maketitle

\begin{center}
$^a$ Faculty of Natural Sciences and Mathematics, University of Maribor, Slovenia\\
\medskip

$^b$ Institute of Mathematics, Physics and Mechanics, Ljubljana, Slovenia\\
\medskip

$^c$ Faculty of Information Technology, University of Pannonia, Veszpr\'em, Hungary\\
\medskip

$^d$ Alfr\'ed R\'enyi Institute of Mathematics, Hungarian Academy of
Sciences, Budapest, Hungary
\medskip

$^e$ Faculty of Mathematics and Physics, University of Ljubljana, Slovenia\\
\medskip

$^f$ Abelium R\&D, Ljubljana, Slovenia\\
\medskip

\end{center}

\maketitle

\begin{abstract}
In a graph $G$ a sequence $v_1,v_2,\dots,v_m$ of vertices is Grundy
dominating if for all $2\le i \le m$ we have $N[v_i]\not\subseteq
\cup_{j=1}^{i-1}N[v_j]$ and is Grundy total dominating if for all
$2\le i \le m$ we have $N(v_i)\not\subseteq \cup_{j=1}^{i-1}N(v_j)$.
The length of the longest Grundy (total) dominating sequence has
been studied by several authors. In this paper we introduce two
similar concepts when the requirement on the neighborhoods is
changed to $N(v_i)\not\subseteq \cup_{j=1}^{i-1}N[v_j]$ or
$N[v_i]\not\subseteq \cup_{j=1}^{i-1}N(v_j)$. In the former case we
establish a strong connection to the zero forcing number of a graph,
while we determine the complexity of the decision problem in the
latter case. We also study the relationships among the four
concepts, and discuss their computational complexities.
\end{abstract}



\section{Introduction}

In~\cite{bgmrr-2014} the Grundy domination number $\ggr(G)$ of a
graph $G$ was introduced as the length of a longest dominating
sequence, that is, a sequence of vertices, such that each vertex of
the sequence dominates at least one new vertex. While the problem is
NP-complete in general, it was demonstrated in the paper
that $\ggr$ can be obtained efficiently within the classes of trees,
split graphs, and cographs. Interval graphs were added to this list
in~\cite{bgk2016+} where in addition a closed formula for the Grundy
domination number of Sierpi\'nski graph was proved. The
investigation of the Grundy domination number on standard graph
products, mostly products of paths and cycles, was done
in~\cite{bbgkkptv2016}. It was also observed that the edge clique
cover number is an upper bound for the Grundy domination number.

It is quite common that along with a domination concept its total
counterpart is studied. And indeed, the Grundy total domination
number $\ggrt(G)$ of a graph $G$ was studied in~\cite{bhr-2016}.
This concept is closely related to the hypergraph of open
neighborhoods of a graph. Using a connection with edge covering
sequences in hypergraphs the authors proved that also the total
Grundy domination number decision problem is NP-complete. Among
several additional results proved we extract that $\ggrt(G) \le
2\ggr(G)$ holds for any graph without isolated vertices.

Motivated by the above two concepts, we introduce in this paper the
closely related concept of the Z-Grundy domination number $\ggrz(G)$
of a graph $G$. We are going to prove in
Section~\ref{sec:zero-forcing} that to determine $\ggrz(G)$ is
equivalent to compute $Z(G)$, where $Z(G)$ is the extensively
studied zero forcing number~\cite{AIM} which is in turn closely
related to the concept of power domination~\cite{BFF-2016+,
haynes-2002}. The complexity of the decision problem whether $Z(G)$ is at least some constant was shown to be NP-complete by~\cite{Aazami2008} and the zero forcing
number has been determined for many classes of graphs~\cite{AIM,barioli, huang, taklimi}.
 This connection on the one hand enables us to deduce
$\ggrz(G)$ for many classes of graphs, and to determine the zero
forcing number of Sierpi\'nski graphs and lexicographic products of paths and cycles on the other hand. To make the
picture complete we introduce in Section~\ref{sec:L-grundy} the
L-Grundy domination number, and prove relations between the four
types of the Grundy domination number. In Section~\ref{sec:complexity-L-grundy} we prove
that the decision problem for the L-Grundy domination number is
NP-complete even when restricted to bipartite graphs. To prove this
result hypergraphs again turn out to be useful.

In the rest of the introduction we recall the definitions of the
Grundy domination number and the Grundy total domination number, and
list some related terminology.

Let $S=(v_1,\ldots,v_k)$ be a sequence of distinct vertices of a graph $G$. The corresponding set $\{v_1,\ldots,v_k\}$ of vertices from the sequence $S$ will be denoted by $\widehat{S}$. The initial segment $(v_1,\dots,v_i)$ of $S$ will be denoted by $S_i$. A sequence $S=(v_1,\ldots,v_k)$, where $v_i\in V(G)$, is called a {\em (legal) closed neighborhood sequence} if, for each $i$
\begin{equation}
\label{eq:defGrundy}
N[v_i] \setminus \bigcup_{j=1}^{i-1}N[v_j]\not=\emptyset.
\end{equation}
We also say that $v_i$ is a {\em legal choice for Grundy domination}, when the above inequality holds.
If for a closed neighborhood sequence $S$, the set $\widehat{S}$ is a dominating set of $G$, then $S$ is called a {\em dominating sequence} of $G$.
We will use a suggestive term by saying that $v_i$ {\em footprints} the vertices from $N[v_i] \setminus \bigcup_{j=1}^{i-1}N[v_j]$, and that $v_i$ is the {\em footprinter} of any $u\in N[v_i] \setminus \bigcup_{j=1}^{i-1}N[v_j]$. For a dominating sequence $S$ any vertex in $V(G)$ has a unique footprinter in $\widehat{S}$.
Clearly the length $k$ of a dominating sequence $S=(v_1,\ldots,v_k)$ is bounded from below by the domination number $\gamma(G)$ of a graph $G$. We call the maximum length of a legal dominating sequence in $G$ the {\em Grundy domination number}
of a graph $G$ and denote it by $\gamma_{gr}(G)$. The corresponding sequence is called a {\em Grundy dominating sequence} of $G$ or $\gamma_{gr}$-sequence of $G$.

In a similar way total dominating sequences were introduced in~\cite{bhr-2016}, for graphs without isolated vertices.
Using the same notation as in the previous paragraph, we say that
a sequence $S=(v_1,\ldots,v_k)$, where $v_i\in V(G)$, is called a {\em (legal) open neighborhood sequence} if, for each $i$
\begin{equation}
\label{eq:defGrundyT}
N(v_i) \setminus \bigcup_{j=1}^{i-1}N(v_j)\not=\emptyset.
\end{equation}
If the above is true, each $v_i$ is said to be a {\em legal choice for Grundy total domination}, and we speak of {\em total footprinters} of which the meaning should be clear. It is easy to see that an open neighborhood sequence $S$ in $G$ of maximum length yields $\widehat{S}$ to be a total dominating set; the sequence $S$ is then called a {\em Grundy total dominating sequence} or {\em $\ggrt$-sequence}, and the corresponding invariant the {\em Grundy total domination number} of $G$, denoted $\ggrt(G)$. Any legal open neigborhood sequence $S$, where $\widehat{S}$ is a total dominating set is called a {\em total dominating sequence}.

\section{Z-Grundy domination and zero forcing}
\label{sec:zero-forcing}

In view of~\eqref{eq:defGrundy} and~\eqref{eq:defGrundyT} the following definition is natural. Let $G$ be a graph without isolated vertices. A sequence $S=(v_1,\ldots,v_k)$, where $v_i\in V(G)$,
is called a {\em (legal) Z-sequence} if, for each $i$
\begin{equation}
\label{eq:defGrundyZ}
N(v_i) \setminus \bigcup_{j=1}^{i-1}N[v_j]\not=\emptyset.
\end{equation}
The {\it{Z-Grundy domination number}} $\ggrz(G)$ of the graph $G$ is the length of a longest Z-sequence. Note that such a sequence is a legal closed neighborhood sequence and hence $\ggrz(G) \leq \ggr(G)$. Given a Z-sequence $S$, the corresponding set $\widehat{S}$  of vertices will be called a {\it{Z-set}}.
\begin{rem}\label{rem}
For any graph $G$ without isolated vertices the equality $\ggrz(G)=\ggr(G)$ holds if and only if there exists a Grundy dominating sequence for $G$ each vertex of which footprints some of its neighbors.
\end{rem}

In this section we connect the Z-Grundy domination number with the well-established zero forcing number that is defined below. (This is also the reason for using ``Z" in the name of the concept.)

Let $G$ be a simple graph with vertex set $V(G)=\{1,\ldots , n\}.$
The {\em minimum rank} mr($G$) of $G$ is the smallest possible rank
over all symmetric real matrices whose $(i,j)$-th entry, for $i \neq j$,
is nonzero whenever vertices $i$ and $j$ are adjacent in $G$ and is
zero otherwise. The {\em{maximum nullity}} M($G$) of $G$ is the
biggest possible nullity over all the above matrices. Clearly
M($G$)$+$mr($G$)$=|V(G)|$. Next, we present a concept derived from
the minimum rank of a graph, cf.~\cite{AIM}. Suppose that the
vertices of a graph $G$ are colored white and blue. If a given blue
vertex has exactly one white neighbor $w$, then by the {\it{color
change rule}} the color of $w$ is changed to blue. A {\em zero
forcing set} for $G$ is a subset $B$ of its vertices such that if
initially vertices from $B$ are colored blue and the remaining
vertices are colored white, then by a repeated application of the
color change rule all the vertices of $G$ are turned to blue. The
{\it{zero forcing number}} Z$(G)$ of a graph $G$ is the size of a
minimum zero forcing set. The zero forcing number is closely related
to minimum rank (maximum nullity) because $|V(G)|-$mr($G$)$\leq
Z(G)$, see~\cite{AIM}.

The main result of this section reads as follows.

\begin{thm}\label{th:zf vs grz}
If $G$ is a graph without isolated vertices, then
$$\ggrz(G)+Z(G)=|V(G)|\,.$$
Moreover, the complement of a (minimum) zero forcing set of $G$ is a
(maximum) Z-set of $G$ and vice versa.
\end{thm}

\proof
Without loss of generality we may assume that $G$ is connected. Let $n=|V(G)|.$

Let $B$ be a zero forcing set, i.e.\ the set of vertices initially
colored blue, and let $k=n-|B|$. The following two sequences appear
in the color change process in which eventually all vertices become
blue. These are the sequences $(b_1,\ldots , b_k)$ and $(w_1,\ldots
, w_k)$, where $b_i$ is the blue vertex selected in the
$i^{\textrm{th}}$ step of the color change process, and $w_i$ is its
unique white neighbor at that moment. We now claim that $(w_k,\ldots
,w_1)$ is a Z-sequence. Note that $b_k\neq w_k$ and hence $w_k$
footprints $b_k$. Let now $i<k$. Then $N(b_i) \subseteq
V(G)\setminus\{w_{i+1},\ldots , w_k\}$, because otherwise in the
$i^{\textrm{th}}$ step of the color change process the blue vertex
$b_i$ would be adjacent to the white vertex $w_i$ and at least one
more white vertex from  the set $\{w_{i+1},\ldots , w_k\},$ which
contradicts the color change rule. It follows that $w_i$ footprints
$b_i$. In particular, if $B$ is a minimum zero forcing set, that is,
$|B|=n-k=Z(G)$, then $k \leq \ggrz(G)$ which in turn implies that
$\ggrz(G) \geq n-Z(G).$

Conversely, let $S=(u_1,\ldots , u_k)$ be a Z-sequence for $G$. We
claim that $X=V(G)-\widehat{S}$ is a zero forcing set. Let
$(a_1,\ldots , a_k)$ be a sequence, where $a_i$ is an arbitrary
vertex selected from the neighbors of $u_i$ that are footprinted by
$u_i$. At the beginning we color all vertices from $X$ blue and the
remaining vertices white. Note that $a_k$ is colored blue because it
was footprinted by the last vertex of $S$. Since $u_k$ is the only
white neighbor of $a_k$, we can color $u_k$ in the first step of the
color change process blue. Let now $i$ be an arbitrary index smaller
than $k$ and assume that the vertices $u_j$, where $i<j\leq k$, are
already colored blue. Since $a_i \in X \cup \{u_{i+1},\ldots ,
u_k\}$, we infer that at this moment $a_i$ is colored blue. Now, by
the same argument as we used for $u_k$, we color $u_i$ blue. By this
procedure we end up with all vertices colored blue, which implies
that $X$ is a zero forcing set. In particular, if $S$ is a maximum
Z-sequence, that is $|\widehat{S}|=\ggrz(G)=k$, then since $Z(G)
\leq n-k$ we get $\ggrz(G) \leq n-Z(G).$ \qed

Combining Theorem~\ref{th:zf vs grz} with the inequality $\ggrz(G) \leq \ggr(G)\,$ we get the following.

\begin{cor}
\label{cr:gr vs zf}
If $G$ is a graph without isolated vertices, then $$\ggr(G)+Z(G) \geq |V(G)|\,.$$
In particular, if $G$ has a Grundy dominating sequence that is also a Z-sequence, then $$\ggr(G)+Z(G)=|V(G)|\,.$$
\end{cor}

\proof
The first assertion follows because $\ggrz(G) \leq \ggr(G)$ holds for any graph $G$. If in addition $S$ is a Grundy dominating sequence that is also a Z-sequence, then $\ggr(G)=|\widehat{S}| \leq \ggrz(G)$ and hence $\ggr(G)=\ggrz(G).$
\qed

There are many examples where equality holds in Corollary~\ref{cr:gr vs zf}. For instance, $\ggr(P_n) = n-1 =\ggrz(P_n)$ and $Z(P_n)=1.$
To see that  $\ggr(G)+Z(G)$  can be arbitrary larger than the order of $G$, consider the stars $K_{1,n}$. Note that $\ggr(K_{1,n})=n$ and that $Z(K_{1,n})=n-1,$ so that $\ggr(K_{1,n})+Z(K_{1,n})=2n-1.$ On the other hand, $\ggrz(K_{1,n})=2$ in accordance with Theorem~\ref{th:zf vs grz}.

The above results are two-fold useful. If $Z(G)$ is known, then
Theorem~\ref{th:zf vs grz} gives us $\ggrz(G)$ and a lower on $\ggr(G)$.  On the other hand,
since the upper bound on the zero forcing number is usually obtained
by construction, there are more problems with lower bounds. But if
$\ggr(G)$ is known for a graph $G$,  then $|V(G)|-\ggr(G)$ is a
lower bound for $Z(G)$. Moreover, if there exists a Grundy
dominating sequence that is a Z-sequence, then the value of $Z(G)$
immediately follows by Corollary~\ref{cr:gr vs zf}. We continue with
some applications of these two approaches.

The Cartesian product $G\cp H$, the strong product $G\boxtimes H$, and the lexicographic prodict $G \circ H$ of
graphs $G$ and $H$ all have the vertex set $V(G)\times V(G)$. In $G\cp H$ vertices $(g,h)$ and $(g',h')$ are adjacent when ($gg'\in E(G)$ and $h=h'$) or ($g=g'$ and $hh'\in E(H)$). The strong product $G\boxtimes H$ is obtained from $G\cp H$ by adding all the edges of the form $(g,h)(g',h')$, where $gg'\in E(G)$ and $hh'\in E(H)$; cf.~\cite{hik-2011}. In the lexicographic product $G\circ H$, vertices
$(g,h)$ and $(g',h')$ are adjacent if either $gg'\in E(G)$ or ($g=g'$ and $hh'\in E(H)$).

The formulas for the zero forcing numbers of some graph products
derived in~\cite{AIM} (for the products in items (i)--(vii) below),
and in~\cite{BFF-2016+} (for item (viii)), together with
Theorem~\ref{th:zf vs grz}, imply the following formulas for their
Z-Grundy domination numbers.

\begin{cor}
\label{cor:many-exact}
\begin{enumerate}[(i)]
\item If $n\ge 1$, then $\ggrz(Q_n) = 2^{n-1}$.
\item If $s, t\ge 2$, then $\ggrz(K_s\cp P_t) = s(t - 1)$.
\item If $2\le |V(G)|\le t$, then $\ggrz(G\cp P_t) = |V(G)|(t-1)$.
\item If $2\le s\le t$, then $\ggrz(P_s\cp P_t) = s(t - 1)$.
\item If $s\ge 3$ and $t\ge 2$, then $\ggrz(C_s\cp P_t) = st - \min\{s,2t\}$.
\item If $s\ge 4$ and $t\ge 2$, then $\ggrz(C_s\cp K_t) = t(s-2)$.
\item If $s, t\ge 2$, then $\ggrz(P_s\boxtimes P_t) = (s-1)(t-1)$.
\item If $3\le s\le t$ and $(s,t)\neq (2r+1,2r+1)$ for some $r \ge 1$, then $\ggrz(C_s\cp C_t)=st-2s$. Moreover, $\ggrz(C_s\cp C_s)=s^2-2s+1$ if $s$ is odd.
\end{enumerate}
\end{cor}

Formulas in items (iv), (v),(vii), and (viii) were independently established
in~\cite{bbgkkptv2016}. More precisely, these formulas were obtained
for the Grundy domination number, yet the sequences applied were
also Z-sequences. Hence the results follow by Remark \ref{rem}.

Next, we add the following results on the zero forcing number of
some product graphs for which the Z-Grundy domination number was
established earlier.

\begin{prp}
\label{prop:cycle-cycle}
\begin{enumerate}[(i)]
\item
If $s\ge 3$ and $t\ge 2$, then $Z(C_s\boxtimes P_t) = 2t + s -2$.
\item
If $s,t > 2$, then
\begin{displaymath}
Z(P_s \circ P_t)=
\left\{ \begin{array}{l l}
\frac{s}{2}\cdot (t+1)-1, & \textrm{$s$ is even}\\
st-\left\lceil \frac{s}{2} \right\rceil \cdot (t-1), & \textrm{$s$ is odd.}\\
\end{array}
\right.
\end{displaymath}
\item
If $s,t > 2$, then
\begin{displaymath}
Z(P_s \circ C_t)=
\left\{ \begin{array}{l l}
\frac{s}{2}\cdot (t+2)-1, & \textrm{$s$ is even}\\
st-\left\lceil \frac{s}{2} \right\rceil \cdot (t-2), & \textrm{$s$ is odd.}\\
\end{array}
\right.
\end{displaymath}
\item
If $s,t > 3$, then
\begin{displaymath}
Z(C_s \circ C_t)=
\left\{ \begin{array}{l l}
\frac{s}{2}\cdot (t+2), & \textrm{$s$ is even}\\
st-\left\lfloor  \frac{s}{2} \right\rfloor \cdot (t-2) -1, & \textrm{$s$ is odd.}\\
\end{array}
\right.
\end{displaymath}
\end{enumerate}
\end{prp}

\proof
These follow by a similar argument to the above. Grundy domination number of $C_s\boxtimes P_t$, $P_s \circ P_t$, $P_s \circ C_t$, and $C_s \circ C_t$ were determined in ~\cite[Corollary 28, Corollary 11, Corollary 12, Corollary 14]{bbgkkptv2016} and all lower bounds use Grundy sequences that are Z-Grundy sequences as well. Therefore Corollary~\ref{cr:gr vs zf} finishes the proof in all cases.
\qed

For another application consider the Sierpi\'nski graphs introduced in~\cite{klavzar-1997} and extensively studied by now, cf.\ the recent survey~\cite{hinz-2017}. They are defined as follows. Set $[n] = \{1,2,\dots,n \}$ and $[n]_0 = \{0,1,\dots,n-1 \}$. The Sierpi\'nski graph $S_p^n$ ($n,p \ge 1$) is defined on the vertex set $[p]_0^n$, two different vertices $u = (u_1,u_2,\dots,u_n)$ and $v = (v_1,v_2,\dots,v_n)$ being adjacent if and only if there exists an $h \in [n]$ such that
\begin{enumerate}
    \item $u_t = v_t$, for $t = 1,2,\dots,h-1$;
    \item $u_h \neq v_h$; and
    \item $u_t = v_h$ and $u_h = v_t$ for $t = h+1,h+2,\dots,n$;
\end{enumerate}

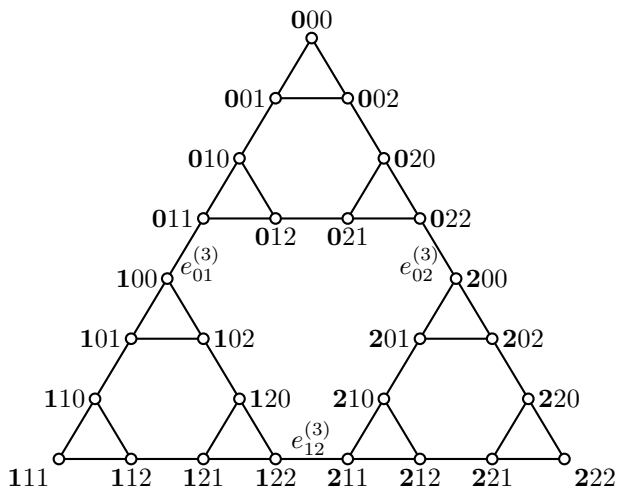
\begin{figure}[ht!]
    \begin{center}
        \begin{tikzpicture}[scale=0.8,style=thick,x=1cm,y=1cm]
        \def\vr{2.5pt} 
        \draw (8,9)--(7.4,8);
        \draw (9.2,5)--(10.4,5);
        \draw (11.6,9)--(12.2,8);
        \draw (9.8,12)--(9.2,11)--(10.4,11)--cycle;
        \draw (8.6,10)--(8,9)--(9.2,9)--cycle;
        \draw (11,10)--(10.4,9)--(11.6,9)--cycle;
        \draw (9.2,11)--(8.6,10);
        \draw (9.2,9)--(10.4,9);
        \draw (10.4,11)--(11,10);
        \draw (9.8,12) [fill=white] circle (\vr);
        \draw (9.2,11) [fill=white] circle (\vr);
        \draw (10.4,11) [fill=white] circle (\vr);
        \draw (8.6,10) [fill=white] circle (\vr);
        \draw (8,9) [fill=white] circle (\vr);
        \draw (9.2,9) [fill=white] circle (\vr);
        \draw (11,10) [fill=white] circle (\vr);
        \draw (10.4,9) [fill=white] circle (\vr);
        \draw (11.6,9) [fill=white] circle (\vr);
        \draw[anchor = south] (9.8,12) node {{\bf 0}00};
        \draw[anchor = east] (9.2,11) node {{\bf 0}01};
        \draw[anchor = west] (10.4,11) node {{\bf 0}02};
        \draw[anchor = east] (8.6,10) node {{\bf 0}10};
        \draw[anchor = east] (8,9) node {{\bf 0}11};
        \draw[anchor = north] (9.2,9) node {{\bf 0}12};
        \draw[anchor = west] (11,10) node {{\bf 0}20};
        \draw[anchor = north] (10.4,9) node {{\bf 0}21};
        \draw[anchor = west] (11.6,9) node {{\bf 0}22};
        \draw (7.4,8)--(6.8,7)--(8,7)--cycle;
        \draw (6.2,6)--(5.6,5)--(6.8,5)--cycle;
        \draw (8.6,6)--(8,5)--(9.2,5)--cycle;
        \draw (6.8,7)--(6.2,6);
        \draw (6.8,5)--(8,5);
        \draw (8,7)--(8.6,6);
        \draw (7.4,8) [fill=white] circle (\vr);
        \draw (6.8,7) [fill=white] circle (\vr);
        \draw (8,7) [fill=white] circle (\vr);
        \draw (6.2,6) [fill=white] circle (\vr);
        \draw (5.6,5) [fill=white] circle (\vr);
        \draw (6.8,5) [fill=white] circle (\vr);
        \draw (8.6,6) [fill=white] circle (\vr);
        \draw (8,5) [fill=white] circle (\vr);
        \draw (9.2,5) [fill=white] circle (\vr);
        \draw[anchor = east] (7.4,8) node {{\bf 1}00};
        \draw[anchor = east] (6.8,7) node {{\bf 1}01};
        \draw[anchor = west] (8,7) node {{\bf 1}02};
        \draw[anchor = east] (6.2,6) node {{\bf 1}10};
        \draw[anchor = north east] (5.6,5) node {{\bf 1}11};
        \draw[anchor = north] (6.8,5) node {{\bf 1}12};
        \draw[anchor = west] (8.6,6) node {{\bf 1}20};
        \draw[anchor = north] (8,5) node {{\bf 1}21};
        \draw[anchor = north] (9.2,5) node {{\bf 1}22};
        \draw (12.2,8)--(11.6,7)--(12.8,7)--cycle;
        \draw (11,6)--(10.4,5)--(11.6,5)--cycle;
        \draw (13.4,6)--(12.8,5)--(14,5)--cycle;
        \draw (11.6,7)--(11,6);
        \draw (11.6,5)--(12.8,5);
        \draw (12.8,7)--(13.4,6);
        \draw (12.2,8) [fill=white] circle (\vr);
        \draw (11.6,7) [fill=white] circle (\vr);
        \draw (12.8,7) [fill=white] circle (\vr);
        \draw (11,6) [fill=white] circle (\vr);
        \draw (10.4,5) [fill=white] circle (\vr);
        \draw (11.6,5) [fill=white] circle (\vr);
        \draw (13.4,6) [fill=white] circle (\vr);
        \draw (12.8,5) [fill=white] circle (\vr);
        \draw (14,5) [fill=white] circle (\vr);
        \draw[anchor = west] (12.2,8) node {{\bf 2}00};
        \draw[anchor = east] (11.6,7) node {{\bf 2}01};
        \draw[anchor = west] (12.8,7) node {{\bf 2}02};
        \draw[anchor = east] (11,6) node {{\bf 2}10};
        \draw[anchor = north] (10.4,5) node {{\bf 2}11};
        \draw[anchor = north] (11.6,5) node {{\bf 2}12};
        \draw[anchor = west] (13.4,6) node {{\bf 2}20};
        \draw[anchor = north] (12.8,5) node {{\bf 2}21};
        \draw[anchor = north west] (14,5) node {{\bf 2}22};
        \draw (7.95,8.25) node {$e_{01}^{(3)}$};
        \draw (11.60,8.25) node {$e_{02}^{(3)}$};
        \draw (9.8,5.35) node {$e_{12}^{(3)}$};
        \end{tikzpicture}
    \end{center}
    \caption{The Sierpi\'nski graph $S_3^{3}$}
    \label{fig:S3}
\end{figure}

It quickly follows from the definition that $S_p^n$ can be constructed from $p$ copies of $S_p^{n-1}$ as follows. For each $i\in [p]_0$ concatenate $i$ to the left of the vertices in a copy of $S_p^{n-1}$ and denote the obtained graph with $iS_p^{n-1}$. Then for each $i\neq j$ join copies $iS_p^{n-1}$ and $jS_p^{n-1}$ by the single edge $e_{ij}^{(n)} = \{ij^{n-1}, ji^{n-1}\}$. $S_3^3$ is depicted in Figure~\ref{fig:S3}.

\begin{prp}
If $n \geq 1$ and $p \geq2$, then $$Z(S_p^n)=\frac{p}{2}\left(p^{n-2}(p-2)+1\right)\,.$$
\end{prp}

\proof
It was proved in~\cite[Theorem 3.1]{bgk2016+} that
\begin{equation*}\label{eq:Sierp}
\ggr(S_p^n)=p^{n-1}+\frac{p(p^{n-1}-1)}{2}\,.
\end{equation*}
In the proof of this theorem it is shown that every vertex added to the appropriate Grundy dominating sequence footprints one of its neighbors. Consequently, the constructed Grundy sequence is also a Z-sequence. Therefore by Corollary~\ref{cr:gr vs zf}, $\ggr(S_p^n)=\ggrz(S_p^n),$ and hence
\begin{eqnarray*}
Z(S_p^n)&=&|V(S_p^n)|-\ggr(S_p^n)= p^n-\left(p^{n-1}+\frac{p(p^{n-1}-1)}{2}\right)\\
&=&\frac{p}{2}\left(p^{n-2}(p-2)+1\right)\,.
\end{eqnarray*}
\qed

\section{L-Grundy domination number and relations between the four concepts}
\label{sec:L-grundy}

Observing the defining conditions in~\eqref{eq:defGrundy}-\eqref{eq:defGrundyZ} it appears natural to consider the remaining, forth related concept. It gives the longest sequences among all four versions, and we call it L-Grundy domination. Given a graph $G$, a sequence $S=(v_1,\ldots,v_k)$ of {\bf distinct} vertices from $G$ is called a {\em (legal) L-sequence} if, for each $i$
\begin{equation}
\label{eq:defGrundyL}
N[v_i] \setminus \bigcup_{j=1}^{i-1}N(v_j)\not=\emptyset.
\end{equation}
Then the {\it{L-Grundy domination number}}, $\ggrl(G)$, of the graph $G$ is the length of a longest L-sequence. Given an L-sequence $S$, the corresponding set $\widehat{S}$  of vertices will be called an {\it{L-set}} (the requirement that all vertices in $S$ are distinct prevents the creation of an infinite sequence by repetition of one and the same vertex). Note that it is possible that some $v_i$ in an L-sequence L-footprints only itself, but at the same time it does not totally dominate any vertex. For instance, consider the star $K_{1,n}$, $n\ge 2$, with the center $w$ and leaves $v_1, \ldots, v_n$. Then $(v_1,\ldots, v_n,w)$ is an L-sequence in which each of the vertices $v_2, \ldots, v_n$ L-footprints only itself. Note also that if a vertex footprints itself, then it will be footprinted later by some other vertex.

The basic inequalities about the four concepts are collected in the
next result.

\begin{prp}
\label{prp:sharpness}
If $G$ is a graph with no isolated vertices, then
\begin{enumerate}[(i)]
\item $\ggrz(G)\le\ggr(G) \le \ggrl(G)-1$,
\item $\ggrz(G)\le \ggrt(G) \le \ggrl(G)$,
\end{enumerate}
and all the bounds are sharp.
\end{prp}

\proof
The bounds $\ggrz(G)\le\ggr(G)$, $\ggrz(G)\le \ggrt(G)$, and $\ggrt(G) \le \ggrl(G)$ follow directly from definitions. To prove that $\ggr(G)\le \ggrl(G)-1$ holds for any graph $G$ with no isolated vertices, consider an arbitrary Grundy dominating sequence $S$ of a graph $G$, and let $v_k$ be the last vertex in the sequence. Then $v_k$ or some of its neighbors has not yet been dominated in previous steps; let us denote the vertex not dominated before the last step by $v_r$, and let $S'$ be the sequence obtained from $S$ by deleting the last vertex $v_k$ from it. Note that $S'\oplus (v_r,v_t)$, where $v_t$ is any neighbor of $v_r$, is an L-sequence of length $k+1$. This readily implies the bound $\ggrl(G)\ge \ggr(G)+1$.

The sharpness of $\ggrz(G)\le\ggr(G)$ is demonstrated by Sierpi\' nski graphs $S_p^n$. To see that $\ggrz(G)\le \ggrt(G)$ is sharp consider for instance stars $K_{1,n}$. To establish the sharpness of $\ggrt(G) \le \ggrl(G)$, use \cite[Theorem 4.2]{bhr-2016}, where the graphs $G$ for which $\ggrt(G)=|V(G)|$ holds are characterized.  (A simple family in this class of graphs are paths of even order.) Finally, to show the sharpness of $\ggr(G)\le \ggrl(G)-1$, the class of graphs $G$ for which $\ggr(G)=|V(G)|-1$, considered briefly in \cite{bgmrr-2014}, presents a family with sharpness of this bound (examples are stars $K_{1,n}$ and paths $P_n$).
\qed

\begin{figure}[h!]
\begin{center}
\begin{tikzpicture}[]
\tikzstyle{vertex}=[circle, draw, inner sep=0pt, minimum size=2pt]
\tikzset{vertexStyle/.append style={rectangle}}
    \vertex (1) at (0,0) [label=below:$\ggrz$] {};
    \vertex (2) at (-1,1) [label=left:$\ggrt$] {};
    \vertex (3) at (1,1) [label=right:$\ggr$] {};
    \vertex (4) at (0,2) [label=above:$\ggrl$] {};
        \path
        (1) edge[dashed] (2)
        (1) edge[dashed] (3)
        (2) edge[dashed] (4)
        (3) edge[dashed] (4);
        \draw (0.7,1.6) node {$\mathbf{+1}$}
            ;
\end{tikzpicture}
\end{center}
\caption{Relations between Grundy domination numbers. In the Hasse diagram an invariant is below another invariant if and only if this is true for the respective Grundy domination numbers for any graph for which both invariants are well-defined. The label $\mathbf{+1}$ describes that $\ggrl(G)\ge \ggr(G)+1$ for any graph $G$ with no isolated vertices.}
\label{fig:POGrundy}
\end{figure}
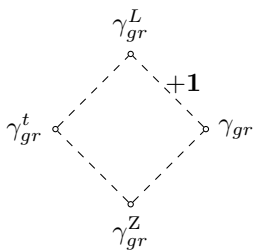

Proposition~\ref{prp:sharpness} can be summarized as the partial order from Figure~\ref{fig:POGrundy}.

As an application of the left inequality in
Proposition~\ref{prp:sharpness}(ii), we get another lower bound for
the zero forcing number of an arbitrary graph with no isolated
vertices. Let $\beta(G)$ stands for the {\em vertex cover
number} of a graph $G$, and $\alpha(G)$ be the {\em independence
number} of $G$.
\begin{prp}\label{prp:beta}
If $G$ is a graph with no isolated vertices, then $Z(G) \ge
\alpha(G)-\beta(G)$.
\end{prp}
\proof Note that Proposition~\ref{prp:sharpness}(ii) implies that
$Z(G)\ge |V(G)|-\ggrt(G)$ for any graph $G$ with no isolated
vertices. By~\cite[Proposition 3.2]{bknt-2016+}, $\ggrt(G)\le
2\beta(G)$, which in turn implies that $Z(G)\ge |V(G)|-2\beta(G)$.
Finally, by the well-known formula $\alpha(G)=|V(G)|-\beta(G)$, we
conclude that $Z(G) \ge \alpha(G)-\beta(G)$.\qed

It is easy to see that the bound $Z(G) \ge \alpha(G)-\beta(G)$ is
sharp in stars $K_{1,n}$ and in odd paths $P_{2k+1}$. It would be
interesting to characterize the extremal graphs for this bound.

It was proved in~\cite{bhr-2016} that $\ggrt(G) \leq 2\ggr(G)$. We can improve this result as follows.

\begin{prp}\label{tz}
If $G$ is a graph without isolated vertices, then $\ggrt(G) \leq 2\ggrz(G)$ and the bound is sharp.
\end{prp}

\proof Let $S=(s_1,\ldots , s_k)$ be a Grundy total dominating sequence and let $A_i$ be the set of vertices totally footprinted by $s_i$. We will prove that at most half of the vertices can be removed from $S$ in such a way that the resulting sequence $S'$ is a legal Z-sequence. If $s_i \in \widehat{S}$ prevents $S$ from being a legal Z-sequence, then $N(s_i)\setminus \bigcup_{j=1}^{i-1}N[s_j] =\emptyset$. Since $S$ is a total dominating sequence this implies $A_i \subseteq  \{s_1,\ldots , s_{i-1}\}$. Let $A=\{s_i \in \widehat{S} : A_i \subseteq \{s_1, \ldots , s_{i-1}\} \}$.

We will first prove that for any $s_i, s_j \in A$, with $i<j$, $A_i \cap A_j = \emptyset$. Suppose that $s_t \in A_i \cap A_j$. Then $t< i <j$. Since $s_t \in A_i$, $s_t$ is totally footprinted by $s_i$ and thus it cannot be totally footprinted by $s_j$, a contradiction.

Finally suppose that for some vertex $s_j \in A$, $s_i \in A \cap A_j$. Then $i < j$. As $s_i \in A$, there exists $t<i$, such that $s_i$ total footprints $s_t$. Therefore, when $s_j$ is added to the sequence, $s_i$ is already totally footprinted by $s_t$, a contradiction with $s_i \in A_j$.

Thus for any $s_i \in A$ there exists $s_t \in (\widehat{S}\setminus A)\cap A_i$ that is not contained in $\bigcup_{j \neq i} A_j$. Therefore $|A| \leq \frac{k}{2}$ and if we define $S'$ to be the sequence obtained from $S$ by removing all its elements that belong to $A$, then $S'$ is a legal Z-sequence, $\ggrz(G) \geq k-|A| \geq k- \frac{k}{2}=\frac{1}{2}\ggrt(G)$.

The bound in Proposition~\ref{tz} is sharp. For example $\ggrt(K_n)=2$ and $\ggrz(K_n)=1$. Another example is a graph $G$ obtained by gluing together two copies of $K_n$ at one vertex. In this case $\ggrt(G)=4=2\ggrz(G)$. \qed

A similar inequality holds between the L-Grundy domination number and the Grundy domination number. Its proof also proceeds along the same lines as the proof of Proposition~\ref{tz}.

\begin{prp}\label{grundylVSgrundy}
If $G$ is a graph, then $\ggrl(G) \leq 2\ggr(G)$ and the bound is sharp.
\end{prp}

\proof
Let $S=(s_1,\ldots , s_k)$ be a legal L-sequence with $\ggrl(G)=k$ and let $A_i$ be the set of vertices L-footprinted by $s_i$. We will prove that at most half of the vertices can be removed from $S$ in such a way that the resulting sequence $S'$ is a legal closed neighborhood sequence. If $s_i \in \widehat{S}$ prevents $S$ from being a legal closed neighborhood sequence, then $N[s_i]\setminus \bigcup_{j=1}^{i-1}N[s_j] =\emptyset$. Since $S$ is an L-sequence this implies that $A_i \subseteq \{s_1,\ldots , s_{i-1}\}$. Let $A=\{s_i \in \widehat{S} : A_i \subseteq \{s_1, \ldots , s_{i-1}\} \}$.

We will first prove that for any $s_i, s_j \in A$, with $i<j$, $A_i \cap A_j = \emptyset$. Suppose that $s_t \in A_i \cap A_j$. Then $t< i < j$. Since $s_t \in A_i$, $s_t$ is L-footprinted by $s_i$ and thus it can not be L-footprinted by $s_j$, a contradiction.

Finally suppose that for some vertex $s_j \in A$, $s_i \in A \cap A_j$. Then $i < j$. As $s_i \in A$, there exists $t<i$, such that $s_i$ L-footprints $s_t$. Therefore, when $s_j$ is added to the sequence, $s_i$ is already L-footprinted by $s_t$, a contradiction with $s_i \in A_j$.

Thus for any $s_i \in A$ there exists $s_t \in (\widehat{S}\setminus A)\cap A_i$ that is not contained in $\bigcup_{j \neq i} A_j$. Therefore $|A| \leq \frac{k}{2}$and if we define $S'$ to be the sequence obtained from $S$ by removing all its elements that belong to $A$, then $S'$ is a legal closed neighborhood sequence in $G$, $\ggr(G) \geq k-|A| \geq k- \frac{k}{2}=\frac{1}{2}\ggrl(G)$. Again, the bound from Proposition~\ref{grundylVSgrundy} is sharp, as $\ggrl(K_n)=2$ and $\ggr(K_n)=1$. Another example is a graph $G$ obtained by gluing together two copies of $K_n$ at one vertex. In this case $\ggrl(G)=4=2\ggr(G)$. \qed

Similar relations between the Grundy domination number and the Z-Grundy domination number, as well as between the L-Grundy domination number and the Grundy total domination number, do not hold. That is, there exist graphs $G$ with $\ggr(G)>C\ggrz(G)$ and graphs $G$ with $\ggrl(G) > C\ggrt(G)$, where $C$ is an arbitrarily chosen positive constant. For example, if $G$ is the star $K_{1,n}$, $n\ge 3$, then $\ggr(K_{1,n})=n$, $\ggrl(K_{1,n})=n+1$, and $\ggrz(K_{1,n})=2$. It was noticed in \cite{bhr-2016} that $\ggr(G)$ can be arbitrarily bigger than $\ggrt(G)$, the stars again form a simple family demonstrating this fact.

\section{On the computational complexity of L-Grundy domination}
\label{sec:complexity-L-grundy}

Recall that NP-completeness results for decision versions of $\ggr$ and $\ggrt$ have been known, see~\cite{bgmrr-2014,bhr-2016}, respectively. Furthermore, it was established in~\cite{Aazami2008} that the zero forcing number yields an NP-complete problem, which combined with the formula $\ggrz(G)+Z(G)=|V(G)|$ from Theorem~\ref{th:zf vs grz} implies that the decision version of $\ggrz$ is also NP-complete. Hence it remains to consider the remaining invariant, $\ggrl$, and the corresponding computational complexity problem:

\begin{center}
\fbox{\parbox{0.85\linewidth}{\noindent
{\sc L-Grundy Domination Number}\\[.8ex]
\begin{tabular*}{.93\textwidth}{rl}
{\em Input:} & A graph $G=(V,E)$, and an integer $k$.\\
{\em Question:} & Is $\ggrl(G)\ge k$?
\end{tabular*}
}}
\end{center}

In the study of the above problem we involve covering sequences in hypergraphs,
in a similar way as in~\cite{bgmrr-2014}, where the connection between dominating sequences in graphs and covering sequences in hypergraphs have been established, see also~\cite{bhr-2016} for another similar application of covering sequences in hypergraphs.

An \emph{edge cover} of a hypergraph $\cH=(X,{\cal E})$ with no isolated vertices,
is a set of hyperedges from ${\cal E}$ that cover all vertices of $X$.  That is, the union of the hyperedges from an edge cover is the ground set $X$. The minimum number of hyperedges in an edge cover of $\cal H$ is called the \emph{(edge) covering number} of $\cal H$ and is denoted by $\rho(\cH)$.
A sequence ${\cal C}=(C_1,\ldots,C_r)$, where $C_i\in \cal E$, is called a \emph{(legal) hyperedge sequence} of $\cal H$, if for any $i$, $1\le i \le r$,  $C_i$ is picked in such a way that it covers some vertex not captured by previous steps;
that is, $C_i\setminus (\cup_{j<i}{C_j})\neq \emptyset$.  In this case ${\cal C}$ is called an \emph{edge covering sequence}, and the maximum length $r$ of an edge covering sequence of $\cal H$  is called the \emph{Grundy covering number} of $\cal H$, and is denoted by $\rho_{gr}(\cH)$. It was shown in~\cite{bgmrr-2014} that determining whether $\rho_{gr}(\cH)$ is bounded from below by a given constant is an NP-complete problem.

\begin{thm}
{\sc L-Grundy Domination Number} is NP-complete, even when restricted to bipartite graphs.
\end{thm}
\proof Let $G$ be an arbitrary bipartite graph, with $(A,B)$ being the bipartition of $V(G)$. Let $G^*$ be the bipartite graph obtained from $G$ by adding the set $I$ of $|B|$ independent vertices to the set $A$,  and connecting each vertex of $I$ to each vertex of $B$ by an edge. Clearly, the bipartition of $V(G^*)$ is $(A\cup I,B)$, and $|A\cup I|=|A|+|B|$.

From the (bipartite) graph $G$ we derive the following hypergraph (which is one of the components of the open neighborhood hypergraph of $G$): $$\cH(G,B)=(A,\{N_G(b)\,| \,b\in B\}).$$
In an analogous way $\cH(G^*,B)$ is defined, and clearly $\rho_{gr}(\cH(G,B))=\rho_{gr}(\cH(G^*,B))$.
The following claim is the crucial step in the proof.
\bigskip

\textbf{Claim} $\ggrl(G^*)=|A|+|B|+\rho_{gr}(\cH(G,B))$.

\medskip

\textbf{Proof (of the claim).}
It is easy to see that $\ggrl(G^*)\ge |A|+|B|+\rho_{gr}(\cH(G,B))$. Indeed, consider the sequence, which first takes all vertices of $A\cup I$, and then emulates the Grundy covering sequence $(C_1,\ldots,C_k)$ of $\cH(G^*,B)$ (or equivalently, of $\cH(G,B)$), by choosing $b_i$ as the $i$th term, where $C_i=N_G(b_i)$. It is clear that the resulting sequence is an L-sequence of $G^*$ with length $|A|+|B|+\rho_{gr}(\cH(G,B))$.

For the reversed inequality, let $S$ be an L-sequence in $G^*$. First note that whenever a vertex $b\in B$ appears in $S$, then at most one vertex from $I$ can appear in $S$ after $b$. Now, consider two possibilities. Suppose that the first vertex $b$ from $B$ that appears in $S$ appears in $S$ before the first vertex from $I$. If all vertices from $B$ appear in $S$ before any vertex from $I$ appeared in $S$, then the length of $S$ is at most $|A|+|B|+1   $, which is clearly not bigger than $|A|+|B|+\rho_{gr}(\cH(G,B))$. Otherwise, the sequence $S$ starts with some vertices from $A\cup B$, then a vertex from $I$ appears, and after that some vertices from $A\cup B$ come in $S$. In this case we also get $|\widehat{S}|\le |A|+|B|+1$ as desired.

Second possibility is that the first vertex from $I$ that appears in $S$ appears in $S$ before the first vertex from $B$. Then all vertices from $B$ are totally dominated at the time the first vertex from $B$ is added to $S$. This implies that $|B\cap \widehat{S}|\le \rho_{gr}(\cH(G,B))$, which in turn implies that $|\widehat{S}|\le |A|+|B|+\rho_{gr}(\cH(G,B))$, which finally implies $\ggrl(G^*)\le|A|+|B|+\rho_{gr}(\cH(G,B))$. {\small \qed}

Since determining the Grundy covering number of a hypergraph (with no isolated vertices) is NP-hard, and  $\cH(G,B)$ can represent an arbitrary hypegraph, we infer from the claim that determining $\ggrl(G^*)$ is also NP-hard for any bipartite graph $G$.
\qed

\section{Concluding remarks}

\begin{enumerate}
\item
Chang et al.\ introduced the concept of $k$-power
domination~\cite{cdm-2012} as a natural generalization of power
domination. An analogous concept of the so-called $k$-forcing was
introduced by Amos et al.~\cite{acd-2015}; the latter concept generalizes
the zero forcing number. (For a connection between the two concepts,
see~\cite{fhk-2017+}.) Given a graph $G$, the {\em $k$-forcing set}
is a subset $S\subseteq V(G)$ ($S$ are the vertices initially
colored blue) such that in the color change process all neighbors of
a blue vertex with at most $k$ white neighbors become blue, and at
the end of this propagation procedure all vertices are blue.
Following our approach from Section~\ref{sec:zero-forcing}, we
present the concept of {\em $k$-Z-Grundy domination number} of a
graph $G$, defined as follows.

Let $G$ be a graph with $\delta(G)\ge k$. A sequence
$S=(v_1,\ldots,v_k)$, where $v_i\in V(G)$, is called a {\em (legal)
$k$-Z-sequence} if, for each $i$ there exists a vertex $u_i\in
N(v_i)$, such that $u_i\in N[v_j]$ holds for less than $k$ vertices
$v_j$, $j \in\{1,\ldots,i-1\}$. (For $k=1$ this definition coincides
with that of the (legal) Z-sequence.) Note that the resulting set
$\widehat{S}$ is a $k$-dominating set, i.e., $\widehat{S}$ has
the property that each vertex outside $\widehat{S}$ has at least $k$ neighbors in
$\widehat{S}$. The corresponding invariant, the $k$-domination
number $\gamma_{ k}(G)$ of a graph $G$, was introduced
in~\cite{FJ}, and studied in several papers, see
e.g.~\cite{alt-2015,HV, CFHV}.
 The {\it{$k$-Z-Grundy domination number}} $\gamma_{gr}^{\textrm{Z},k}(G)$
of the graph $G$ is the length of a longest $k$-Z-sequence. By the
above connection with $k$-domination, we infer that
$\gamma_{gr}^{\textrm{Z},k}(G) \ge \gamma_{ k}(G)$.

\item Propagation time for zero forcing was introduced
in~\cite{hhk-2012} as the minimum over all smallest zero-focring
sets of the number of propagation steps. Being in a given state, a
propagation step consist of coloring all vertices blue that may be
colored blue in the state. Analogous concept for the power
domination (under the name propagation radius) was introduced
in~\cite{dk-2014} and independently in~\cite{liao-2016}.

The interpretation in terms of the Z-Grundy domination of these
concepts is the following. White vertices that will turn blue in the
same propagation step have the property that set of the
corresponding blue neighbors partition into private neighbors of
each of the white neighbors, where private neighborhoods are
considered with respect to the set of white vertices, before the
propagation step. In terms of Z-Grundy dominating sequence $S$ that
is built in the reversed process to zero forcing propagation, this
implies that a subsequence of the corresponding (white) vertices can
be permuted in an arbitrary order. Hence, given a zero forcing set,
the number of propagation steps coincides with the number of
uniquely determined consecutive subsequences of the Z-sequence with
the above permutation property.

This leads to the following result.

\begin{prp}
Given a graph $G$, the zero forcing propagation time of $G$
coincides with the minimum number of consecutive permutable
subsequences over all Z-sequences in $G$.
\end{prp}

\item
There are many equality cases, for which it would be interesting to
characterize the extremal graphs. We list some of them:

    $\bullet$ Characterization of graphs $G$ with $|V(G)|=\gamma^L_{gr}(G)$.

    $\bullet$ Characterization of graphs $G$ with $\gamma^L_{gr}(G)=\gamma^Z_{gr}(G)+1$.

    $\bullet$ Characterization of graphs $G$ with $\gamma^L_{gr}(G)=2\gamma_{gr}(G)$.

    $\bullet$ Characterization of graphs $G$ with $\ggrz(G)=\ggr(G)$.

A similar problem about characterizing graphs for which $\gamma^t_{gr}(G)=2\gamma_{gr}(G)$ holds were posed in~\cite{bhr-2016}. All the above questions can be asked in some specific, interesting
families of graphs, such as trees, $k$-regular graphs, and graph products.
In particular, in our next project, we plan to consider different Grundy domination numbers in
grid-like and toroidal graphs with respect to various graphs products.\\

\item
Recall that the zero forcing number is closely related to minimum
rank (and maximum nullity) via the formula: $|V(G)|-$mr($G$)$\leq
Z(G)$. It would be interesting to establish a
(Grundy)-domination-type concept that would be directly connected to
the minimum rank of a graph.
\end{enumerate}

\section*{Acknowledgements}

The authors acknowledge the project (Combinatorial Problems with an Emphasis on Games, N1-0043) was financially supported by the Slovenian Research Agency. The authors acknowledge the financial support from the Slovenian Research Agency (research core funding No.\ P1-0297). The authors acknowledge that the research was in part financially supported by the Slovenian Research Agency; project grant L7-5554.

Research of Cs.\ Bujt\' as, B.\ Patk\' os, Zs.\ Tuza, and M.\ Vizer
was supported by the National Research, Development and Innovation
Office -- NKFIH under the grant SNN 116095.

Research of B.\ Patk\'os was supported by the J\'anos Bolyai Research Fellowship of the Hungarian Academy of Sciences.


\end{document}